\font\tenrm=pbkl8t\font\sevenrm=pbkl8t at 7.6pt\font\fiverm=pbkl8t at 6pt\rm

\input amssym
\font\tt=t1xtt
\font\sc=pbklc8t

\font\mathtext=pbkl7t \fontdimen2\mathtext=0pt
 \font\mathsubtext=pbkl7t at 7.6pt\fontdimen2\mathsubtext=0pt
 \font\mathsubsubtext=pbkl7t at 6pt\fontdimen2\mathsubsubtext=0pt
\textfont0=\mathtext
 \scriptfont0=\mathsubtext
 \scriptscriptfont0=\mathsubsubtext

\font\mathlet=fncmi\fontdimen2\mathlet=0pt
 \font\mathsublet=fncmi at 7.6pt  \fontdimen2\mathsublet=0pt
 \font\mathsubsublet=fncmi at 6pt \fontdimen2\mathsubsublet=0pt
\skewchar\mathlet='177 \skewchar\mathsublet='177 \skewchar\mathsubsublet='177
\textfont1=\mathlet
 \scriptfont1=\mathsublet
 \scriptscriptfont1=\mathsubsublet

\font\mathsym=futsy
 \font\mathsubsym=futsy at 7.6pt
 \font\mathsubsubsym=futsy at 6pt
\skewchar\mathsym='60 \skewchar\mathsubsym='60 \skewchar\mathsubsubsym='60
\textfont2=\mathsym
 \scriptfont2=\mathsubsym
 \scriptscriptfont2=\mathsubsubsym

\textfont4=\tenrm
 \scriptfont4=\sevenrm
 \scriptscriptfont4=\fiverm

\chardef\ss="FF
\chardef\ae="E6
\chardef\oe="F7
\chardef\o="F8
\chardef\l="AA
\chardef\AE="C6
\chardef\OE="D7
\chardef\O="D8
\chardef\L="8A
\chardef\i="19
\chardef\j="1A
\chardef\texttilde="7E

\catcode`@=11

\def\k#1{\char'246}
\def\`#1{{\accent"00 #1}}
\def\'#1{{\accent"01 #1}}
\def\v#1{{\accent"07 #1}} \let\^^_=\v
\def\u#1{{\accent"08 #1}} \let\^^S=\u
\def\=#1{{\accent"09 #1}}
\def\^#1{{\accent"02 #1}} \let\^^D=\^
\def\.#1{{\accent"0A #1}}

\def\~#1{{\accent"03 #1}}
\def\"#1{{\accent"04 #1}}
\def\t#1{{\edef\next{\the\font}\the\textfont1\accent"7F\next#1}}

\mathcode`\(="4428
\mathcode`\)="5429
\mathcode`\:="343A
\mathcode`\;="643B
\mathcode`\[="445B
\mathcode`\]="545D
\mathcode`\+="242B
\mathchardef\colon="643A
\def\l{\char'252}
\def\L{\char'212}
\mathchardef\Gamma="0100
\mathchardef\Delta="0101
\mathchardef\Theta="0102
\mathchardef\Lambda="0103
\mathchardef\Xi="0104
\mathchardef\Pi="0105
\mathchardef\Sigma="0106
\mathchardef\Upsilon="0107
\mathchardef\Phi="0108
\mathchardef\Psi="0109
\mathchardef\Omega="010A

\input amssym

\overfullrule=0pt
\magnification1200
\nopagenumbers
\def\title{Decomposing Bj\"orner's Matrix}
\font\sevenss=pagkc8t at 7pt
\headline{\everymath{\scriptstyle}\sevenss \ifnum\pageno>1
\ifodd\pageno\title\ \hrulefill\ \the\pageno
\else\the\pageno\ \hrulefill\ \'S.~R.~Gal\fi\fi} 
\centerline{\bf \title}
\smallskip
\centerline{\'Swiatos\l aw R. Gal%
\footnote{$^\star$}{Partially supported by a {\sc mnisw} grant
N N201 541738.}}
\centerline{Wroc\l aw University}
\centerline{\tt http://www.math.uni.wroc.pl/\char'176sgal/papers/bjo.pdf}

\bigskip
{
\narrower\narrower\smallskip\noindent\everymath{\scriptstyle}%
\sevenbf Abstract: \sevenrm  We give an alternative proof of a (former) conjecture
of Bj\"orner stating that the matrix expressing face numbers in terms of g numbers
is totally non-negative.  We briefly discuss the case of simple flag polytopes.
\smallskip}   
\bigskip

\parindent0pt
\parskip=\smallskipamount
\footnote{}{2010 {\it Mathematics Subject Classification: } 52B05.}
\footnote{}{{\it Key phrases: g-vector,  Bj\"orner's Conjecture}.}

\parindent0pt
\parskip=\smallskipamount
\newcount\secno\secno=0
\newcount\thmno\thmno=0

\def\section #1\par{\vskip 0pt plus.1\vsize\penalty-250
\vskip 0pt plus-.1\vsize\bigskip\vskip\parskip
\global\advance\secno by 1\global\thmno=0
{\the\secno.} {\sc #1}\par
\smallskip
}

\def\tag{\the\secno.\the\thmno}
\def\mktag{{\global\advance\thmno by 1}\tag}

Let $f(t)=\sum f_it^i$ denote the f-polynomial of a $d$-dimensional simple polytope.
The coefficients ($f_i$ denotes the number of codimension $i$ faces of the polytope)
are called {\sl face numbers} of the polytope.  It is a consequence 
of Dehn-Somerville relations that $f(t)$ can be written a s a linear combination of 
polynomials $u_i(t)=\sum_{q=i}^{d-i}(1+t)^q$.  The
coefficients of the corresponding expantions are denoted $g_i$ i.e.
$$f(t)=\sum_{i=0}^{\lfloor d/2\rfloor}g_iu_i(t).$$
In particular there exists a matrix $M^{(g)}_d$ such that
$$\pmatrix{f_0\cr f_1\cr\vdots\cr f_d}
=M^{(g)}_d\pmatrix{g_0\cr g_1\cr\vdots\cr g_{\lfloor d/2\rfloor}}.$$
It was conjectured by Bj\"orner that $M^{(g)}_d$ is {\sl totally non-negative}  i.e. all its
minors are non-negative.  This was proven by Bj\"orklund and Engstr\"om [BE].

Bj\"orner already proved [Bj] that two-by-two minors of $M^{(g)}_d$ are non-negative and
used it to refine lower and upper bounds for simple polytopes.

In this note we present another proof which showas that $M^{(g)}_d$ is totally non-negative.

Since face and g-numbers are related as follows
$$\eqalign{\sum f_it^{d-i}
&=\sum_{k=0}^{\lfloor d/2\rfloor} g_k\sum_{j=k}^{d-k}(1+t)^j\cr
&=\sum_{k=0}^{\lfloor d/2\rfloor} g_k\sum_i\sum_{j=k}^{d-k}{j\choose d-i}t^{d-i}\cr
&=\sum_i t^{d-i}\sum_{k=0}^{\lfloor d/2\rfloor} g_k\left({d-k+1\choose d-i+1}-{k\choose d-i+1}\right).\cr}$$
We have
$$f_i=\sum_{k=0}^{\lfloor d/2\rfloor}M^{(g)}_d(i,k)g_k$$
with
$$M^{(g)}_d(i,k)={d-k+1\choose d-i+1}-{k\choose d-i+1}.$$

\section Four matrices.

Let us define four infinite matrices:

$$\eqalign{
A_\pm(i,j)&={j+1\choose i-j}\mp{j\choose i-j-1},\cr
G_+(j,k)&={2k+1\over 2j+1}{k+j\choose 2j},\cr
G_-(j,k)&={k+j+1\choose 2j+1}.\cr
}$$
In all the cases we assume that $i$, $j$, and $k$ run through natural numbers inclu\-ding zero.

\proclaim Proposition.  The matrices $A_+$, $A_-$, $G_+$, and $G_-$ are totally
non-negative.

{\it Proof:}
First notice that
$$A_+(i,j)={j\choose i-j},\qquad A_-(i,j)={i+1\over j+1}{j+1\choose i-j}.$$
It is quite standard that
$${j\choose i-j},\quad{j+1\choose i-j},\quad{k+j\choose 2j},\quad{k+j+1\choose 2j+1}$$
are totally non-negative.  For example ${k+j\choose 2j}$ counts the number of paths from $(-2j,j)$
to $(0,k)$ where only steps $(1,0)$ or $(0,1)$ are allowed. (Cf. [FZ, Lemma 1].)

When $X(i,j)$ is totally non-negative, then $a(i)X(i,j)b(j)$ is totally non-negative provided
$a$ and $b$ are non-negative sequences.  Thus the proof.\hfill$\square$

\section Decomposition of the matrix $M^{(g)}$

\proclaim Theorem.
Let $\varepsilon=(-)^d$ and $n=\lfloor d/2\rfloor$.  One has
$$M^{(g)}_d(i,k)=\sum_{j=0}^{n}
A_\varepsilon(i,j)G_\varepsilon(n-j,n-k).$$

\proclaim Remark.  Since total non-negativeness is preserved under reversing
the order of both indices, and multiplication of such matrices total
non-negativity of\/ $M^{(g)}_d$ follows immediately.  We believe that the above
formula is of independent interest. 

{\it Poof of the Theorem:} It reduces to:

$$
\eqalignno{
{2n-k+1\choose 2n-i+1}-{k\choose 2n-i+1}&=
\sum_{j=k}^{n}{j\choose i-j}{2n-2k+1\over 2n-2j+1}{2n-(k+j)\choose 2n-2j}&(1)\cr
{2n-k+2\choose 2n-i+2}-{k\choose 2n-i+2}&=
\sum_{j=k}^{n}{i+1\over j+1}{j+1\choose i-j}{2n-(k+j)+1\choose 2n-2j+1}&(2)\cr}
$$

Let us provide a proof of the first identity.  The second can be obtained in a very
similar way.

Multiply both sides by $s^i$ and sum over $0\leq i\leq 2n+1$.  Left-hand
side of (1) becomes
$$s^{2n+1}(1+1/s)^{2n-k+1}-s^{2n+1}(1+1/s)^k=(1+s)^{2n+1}\left(s\over1+s\right)^k-s^{2n+1}\left(1+s\over s\right)^k.$$
The right-hand side of (1) becomes
$$\sum_{j=k}^{n}(s(s+1))^j{2n-2k+1\over 2n-2j+1}{2n-(k+j)\choose 2n-2j}.$$

Let
$$F_a(z)=\sum_{j=0}^{a}z^j{2a+1\over 2a-2j+1}{2a-j\choose 2a-2j}.$$
Then left-hand side of (1) equals to $(s(s+1))^kF_{n-k}(s(s+1))$.
It is straightforward to check that
$$z(4z+1)F''_a(z)+a(z-a(4z+1))F'_a(z)+2a(2a+1)F_a(z)=0.$$
If $G_a(s)=F_a(s(s+1))$ then
$$s(s+1)G''_a(s)-2a(1+2s)G'_a(s)+2a(2a+1)G(s)=0.$$
The above equation is solved also by $(s+1)^{2a+1}$ and $s^{2a+1}$.
Since $G_a(s)$ and $(s+1)^{2a+1}-s^{2a+1}$ have the same (non-zero) coefficients at $s$ and $s^{2a}$
we conclude that
$$G_a(s)=(s+1)^{2a+1}-s^{2a+1},$$
and, what follows the left-hand side of (1) is equal to
$$(s(s+1))^kG_{n-k}(s)=\left(s\over1+s\right)^k(s+1)^{2n+1}-\left(1+s\over s\right)^ks^{2n+1}.$$
Which completes the proof.\hfill$\square$

\proclaim Question.  The matrices $G_\varepsilon(n-j,n-k)$ depend on dimension
mildly.  And they provide, applied to g-numbers another (natural) expansion of
of f-polynomial.  Does this expansion is combinatorially/geometrically useful?

\section Flag polytopes and $\gamma$-matrix.

If the polytope in question is {\sl flag} (ie. all pairwise intersecting
families of faces are centered) an interesting expansion of the f-vector is given
in terms of polynomials $v_i(t)=(t+1)^i(t+2)^{d-2i}$, ie.
$$f(t)=\sum \gamma_i v_i(t)$$
(cf [Br,G].)

Unsolved question is whether the coefficients $\gamma_i$ are non-negative in this case
([G, Conjecture 2.1.7]).  As before we can write down the matrix which computes the coefficients
of f-polynomial in terms of $\gamma_i$:
$$\eqalign{\sum f_it^{d-i}
&=\sum_{k=0}^{\lfloor d/2\rfloor}\gamma_k(1+t)^k(2+t)^{d-2k}\cr
&=\sum_{k=0}^{\lfloor d/2\rfloor}\gamma_k\sum_p\sum_q{k\choose p}{d-2k\choose q}2^qt^{k-p+d-2k-q}\cr
&=\sum_i t^{d-i}\sum_{k=0}^{\lfloor d/2\rfloor}\gamma_k\sum_q{k\choose i-k-q}{d-2k\choose q}2^q.\cr}$$
Thus we have
$$\pmatrix{f_0\cr f_1\cr\vdots\cr f_d}=M^{(\gamma)}_d\pmatrix{\gamma_0\cr \gamma_1\cr\vdots\cr \gamma_{\lfloor d/2\rfloor}},$$
with
$$M^{(\gamma)}_d(i,k)=\sum_j{k\choose i-k-j}{d-2k\choose j}2^j.$$

One can also express $g_i$ by $\gamma_j$:
$$g_i=\sum_{0\leq j\leq i}\left({n-2j\choose i-j}-{n-2j\choose i-1-j}\right)\gamma_i
=\sum_{0\leq j\leq i}{n-2j+1\choose i-j}{n-2i+1\over n-2j+1}\gamma_i.$$
Since the matrix with entries ${n-2j+1\choose i-j}$ is easily seen to be
totally non-negative we see that $M^{(\gamma)}_d$ is totally non-negative as
$M^{(g)}_d$ is such.

This could be also proven directly by checking that
$$M^{(\gamma)}_d(i,k)=\sum_{j=0}^nA_\varepsilon(i,j)\Gamma(n-j,n-k),$$
where $n=\lfloor d/2\rfloor$, $\varepsilon=(-)^d$, and
$$\Gamma(j,k)=4^{k-j}{k\choose j}.$$

\proclaim Question.  Can one use non-negativity of (say, two-by-two)
minors of $M^{(\gamma)}_d$ to get a (depending on conjectural non-negativity of $\gamma$
numbers) upper/lower bound for face numbers of flag simple polytopes?

\bigskip{\bf References}\par\smallskip

\parindent=.4in

\item{[Bj]} {\sc A.~Bj\"orner},
{\it A comparison theorem for f-vectors of simplicial polytopes},
Pure Appl.\ Math.\ Q. {\bf 3} (2007), no.~1, part 3, pp. 347--356,

\item{[BE]} {\sc M.~Bj\"orklund \& A.~Engstr\"om},
{\it The g-theorem matrices are totally nonnegative},
J.~Combin.\ Theory Ser.~A, {\bf 116} (2009), no.~3, pp. 730--732,

\item{[Br]} {\sc P.~Br\"and\'en}, {\it Sign-graded posets, unimodality of
$W$W-polynomials and the Charney-Davis conjecture}, Electron.\ J.\ Combin.\
{\bf 11}  (2004/06),  no.~2, Research Paper 9, (electronic),

\item{[FZ]} {\sc S.Fomin \& A.Zelevinsky} {\it Total positivity: tests and
parametrizations}, Math.\ Intelligencer {\bf 22} (2000), no. 1, pp. 23--33, 

\item{[G]} {\sc \'S.\ Gal}, {\it Real root conjecture fails for five- and
higher-dimensional spheres}, Discrete Comput.\ Geom.\ {\bf 34} (2005), no.~2,
pp.~269--284.

\parskip=0pt
\bigskip
\obeylines
Mathematical Institute
University of Wroc{\l}aw
pl. Grunwaldzki 2/4
50-384 Wroc{\l}aw
Poland
\tt sgal@math.uni.wroc.pl

\bye